\documentclass[a4paper,12pt]{amsart}
\usepackage{amssymb}
\author{Ivan V.~Arzhantsev}
\address{
Department of Higher Algebra\\
Faculty of Mechanics and Mathematics\\
Moscow State University\\
119992 Moscow, Russia}
\email{arjantse@mccme.ru}
\urladdr{http://mech.math.msu.su/department/algebra/staff/arzhan.htm}
\author{Olga V. Chuvashova}
\address{
Department of Higher Algebra\\
Faculty of Mechanics and Mathematics\\
Moscow State University\\
119992 Moscow, Russia}
\email{chuvashova@yandex.ru}
\thanks{Supported by RFBR grant MAC 03-01-06252, by CRDF
grant RM1-2543-MO-03 and by the RF President grant MK-1279.2004.1}

\title[Homogeneous spaces of complexity one]{Classification of affine
homogeneous spaces of complexity one}

\date{August 8, 2004}

\keywords{Homogeneous spaces, complexity, rank, isotropy
representation, spherical subgroups}

\subjclass[2000]{Primary 14M17, 14R20, 37J35; Secondary 22F30,
32M10, 53C30}

\newcommand{\kk}{\mathbb{K}}
\newcommand{\C}{\mathbb{C}}
\newcommand{\g}{\mathfrak{g}}
\newcommand{\h}{\mathfrak{h}}
\newcommand{\ppp}{\mathfrak{p}}
\newcommand{\m}{\mathfrak{m}}
\newcommand{\bor}{\mathfrak{b}}
\newcommand{\cen}{\mathfrak{z}}
\newcommand{\s}{\mathfrak{s}}
\newcommand{\lev}{\mathfrak{l}}
\newcommand{\pp}{\mathfrak{p}}
\newcommand{\z}{\mathfrak{z}}
\newcommand{\e}{\varepsilon}
\newcommand{\rk}{\mathop{\mathrm{rk}}}
\newcommand{\codim}{\mathop{\mathrm{codim}}}
\newcommand{\p}{\pi}
\newcommand{\ssl}{\mathfrak{sl}}
\newcommand{\ssp}{\mathfrak{sp}}
\newcommand{\so}{\mathfrak{so}}
\newcommand{\ssc}{\mathfrak{c}}
\newcommand{\gl}{\mathfrak{gl}}
\newcommand{\cc}{\mathfrak{c}}
\newcommand{\vv}{\mathfrak{v}}
\newcommand{\ww}{\mathfrak{w}}
\newcommand{\spin}{\mathfrak{spin}}

\newtheorem{theorem}{Theorem}

\newtheorem{proposition}{Proposition}
\newtheorem{lemma}{Lemma}
\theoremstyle{definition}

\theoremstyle{remark}

\begin{document}

\begin{abstract}
The complexity of an action of a reductive algebraic group $G$ on
an algebraic variety $X$ is the codimension of a generic $B$-orbit
in $X$, where $B\subset G$ is a Borel subgroup. We classify affine
homogeneous spaces $G/H$ of complexity one. These results are the
natural continuation of the classification of spherical affine
homogeneous spaces, i.e., spaces of complexity zero.
\end{abstract}

\maketitle

\section{Introduction}

Let $G$ be a connected reductive algebraic group over an
algebraically closed field $\kk$ of characteristic zero and $H$ be
an algebraic subgroup of $G$. The complexity $c(G/H)$ is an
important integer-valued invariant of the homogeneous space $G/H$.
The notion of complexity was introduced by Luna and Vust for
homogeneous spaces~\cite{lv} and by Vinberg for arbitrary
$G$-varieties~\cite{vin1}. It has been shown in further researches
that the complexity plays the key role in the study of
geometry of a homogeneous space, in the theory of its
equivariant embeddings, in the theory of invariant Hamiltonian
systems on the cotangent bundle $T^*(G/H)$, and in other fields of
mathematics concerned with homogeneous spaces.

Let us recall the definition of complexity. Let a reductive group
$G$ act on an irreducible algebraic variety $X$ and  $B$ be a
Borel subgroup of $G$. The \emph{complexity} $c_G(X)$ of the
$G$-variety $X$ is the minimal codimension of $B$-orbits in $X$
for the induced action $B:X$. (We shall write $c(X)$, if it
does not lead to misunderstanding.) By Rosenlicht's Theorem,
$c_G(X)$ equals the transcendence degree of the field $\kk(X)^B$
of rational $B$-invariant functions on~$X$. In addition to
complexity, there is one more important characteristic of the
action. The {\it rank} $\rk(X)$  of an irreducible $G$-variety
$X$ is the rank of the weight lattice of $B$-semiinvariant
rational functions on $X$.

A normal $G$-variety $X$ is called \emph{spherical} if $c(X)=0$
or, equivalently, $\kk(X)^B=\kk$. A space $G/H$ and a subgroup
$H\subseteq G$ are called \emph{spherical} if $G/H$ is a spherical
$G$-variety. The theory of spherical varieties is one of the most
developed parts of the theory of algebraic transformation groups,
see e.g.~\cite{br2}. In particular, for spherical homogeneous
spaces there exists a remarkable theory of
embeddings~\cite{lv},~\cite{knop}, generalizing the theory of
toric varieties. In~\cite{tim}, D.A.~Timashev obtained a similar
(but much more difficult) description of embeddings for
homogeneous spaces of complexity one.

Spaces of small complexity appear naturally in the study of the
question "whether the number of orbits in any embedding of a given
homogeneous space is finite ?". For example, the number of orbits in
any embedding of a spherical homogeneous space is finite . At the
same time, any homogeneous space of complexity $\ge 1$ admits a
projective embedding with infinitely many orbits~\cite{akh1}. More
generally, the maximal value of $G$-modality (that is the maximal
number of parameters in a continuous family of $G$-orbits) over
all embeddings of the space $G/H$ is equal to
$c(G/H)$~\cite{vin1}, \cite{akh2}.

Suppose that the homogeneous space $G/H$ is affine and we consider
only its affine embeddings. It is known that every affine
embedding of the space $G/H$ of complexity $\ge 2$ contains
finitely many $G$-orbits if and only if $G/H$ is affinely closed,
that is any affine embedding of  $G/H$ consists of only one
orbit~\cite{at}. An affine homogeneous space of complexity one
admits an affine embedding with infinitely many orbits if and only
if there exists an extension of $H$ by a one-dimensional torus in
$N_G(H)/H$ such that the complexity of the obtained homogeneous
space is equal to one~\cite{at}. Among homogeneous spaces of
simple groups all the spaces with this property correspond to item
3 of  Table 2 of this paper. In the case of semisimple groups
there are much more such spaces. They correspond to items 9, 10,
15, 17, and 26 of  Table 4. For more information on relations
between complexity and modality, see~\cite{ar}.

Another reasons in favor of the investigation of homogeneous
spaces of small complexity are problems of symplectic geometry and
theory of integrability of invariant Hamiltonian systems. Let $G$
be a connected reductive real Lie group, and $K$ be its connected
closed reductive subgroup. The cotangent bundle $T^*(G/K)$ has a
natural structure of a symplectic variety, and the action $G:G/K$
induces the symplectic action $G:T^*(G/K)$. Denote by $P$ the
momentum map $P: T^*(G/K) \to {\g}^*$ corresponding to this
action. The algebras of real analytic functions $C(T^*(G/K))$ and
$C({\g}^*)$ have the canonical Poisson brackets. Note that for all
$h_1, h_2\in C({\g}^*)$ we have $\{ h_1\circ P, h_2\circ P\}=\{
h_1, h_2\}\circ P$. Functions on $T^*(G/K)$ having the form $h\circ
P$ are called {\it collective}. By Noether's Theorem, collective
functions are integrals for every flow on $T^*(G/K)$ with
$G$-invariant Hamiltonian. In other words, $\{ h\circ P, f\}=0$
for every $f\in C(T^*(G/K))^G$. A $G$-invariant Hamiltonian system
$\dot x=sgrad\,f$ is said to be {\it integrable in the class of
Noether's integrals} or {\it collectively completely integrable},
if there exists a set $\{ h_i\}$, $i=1,\dots,\dim G/K$, of real
analytical functions on $\g^*$ such that the functions $h_i\circ
P$ are functionally independent and every pair is in involution.
It is known that any $G$-invariant Hamiltonian system on
$T^*(G/K)$ is integrable in the class of Noether's integrals if
and only if the complexified homogeneous space $G^{\C}/K^{\C}$ is
spherical~\cite{mik1}. If the subgroup $K$ is compact, then these
conditions are equivalent to the fact that the Poisson
bracket vanishes on every pair of functions from
$C(T^*(G/K))^G$~\cite{gs}, see also ~\cite[Prop.~9]{vin2}.

The recent works ~\cite{MS}, \cite{mik2} are concerned with
investigations of properties of Hamiltonian systems on cotangent
bundles of homogeneous spaces of positive complexity. As shown in
these papers, if $K$ is compact, then the number of independent
$G$-invariant real analytical functions on $T^*(G/K)$, that are
also independent from collective functions, is equal to
$2c(G^{\C}/K^{\C})$. (This number equals the corank of the
symplectic action $G^{\C}:T^*(G^{\C}/K^{\C})$, see~\cite{vin2}).
On the other hand, the maximal number of independent collective
functions in involution equals $\dim G/K - c(G^{\C}/K^{\C})$. In
the case of homogeneous spaces of complexity one
\footnote{In~\cite{mik2} they are called {\it almost spherical
spaces}.}, if we add to a maximal independent system of collective
functions in involution an independent of them invariant function
$F$, then we shall obtain a complete system of functions in
involution. Choosing $F$ commuting with the Hamiltonian, one can
show that any $G$-invariant Hamiltonian system on the cotangent
bundle of a homogeneous space of complexity one is
integrable~\cite[Prop.~12]{mik2}.

We may conclude that homogeneous spaces of small complexity are natural and
important objects. Hence the list of all such spaces $G/H$
might be useful. The case of a reductive subgroup $H$ is a natural
beginning of the classification. (By Matsushima's Criterion, this
condition is equivalent to affinity of the space $G/H$.) Under
this restriction, the list of homogeneous spaces of the complexity
$\le 1$ is not too long. This confirms the hypothesis stated at
the end of~\cite{vin1}, where the problem of classification of
homogeneous spaces of complexity one was posed for the first time.

It is well known that the complexity of a homogeneous space $G/H$
is determined by the pair of the tangent algebras $(\g, \h)$
(Proposition~1). This allows to reduce the classification to the
enumeration of pairs $(\g, \h)$, where $\g$ is a reductive Lie
algebra, $\h$ is its reductive subalgebra such that the generic
codimension of the subspace $\h+\mathrm{Ad}(g)\bor, g\in G$ in the
algebra $\g$ equals one. Without loss of generality, it can be
assumed that the algebra $\g$ is semisimple (Proposition~2).

The aim of this paper is to classify all affine homogeneous spaces
of complexity one for reductive groups in terms of their tangent
algebras. The results are represented in Table 4. The
classification of spherical reductive subgroups in simple groups
was obtained in~\cite{Kr}, the classification of spherical
reductive subgroups in reductive groups was obtained
in~\cite{mik1} and, independently, in~\cite{Brion}, and the
classification of reductive subgroups of complexity one in simple
groups is given in \cite{Pan2}, see also~\cite{MS}. For
convenience of the reader the results of these classifications are
included in the text (see Table~1, Table~3 and Table~2
respectively). Note that the classification is made up to
automorphisms of the ambient algebra $\g$. For example, there
exist exterior automorphisms of $\so_8$ that take the
spherical pairs $(\so_8, \spin_7)$ and $(\so_8,
\ssp_4\oplus\ssl_2)$ to the pairs corresponding to item 10 of
Table~1 with $(n,m)=(7,1)$ ¨ $(5,3)$ respectively, so these pairs
are not represented in the table. In Tables 1 and 2 we also
indicate $\h$-module $\m$ corresponding to the isotropy
representation. In Tables 3 and 4 we  do not indicate the isotropy
representation, but it can be found easily since the pairs of
these tables are obtained by the operation of coupling in the
pairs from Tables 1 and 2. (The definition of this operation will
be given below.) Suppose we make a coupling of a component $\h$;
after this the isotropy representation will be the direct sum of
the isotropy representations of the coupled pairs and one copy of
the adjoint representation of $\h$. In Tables 3 and 4 we indicate
the stationary subalgebras of general position for the isotropy
representation and the rank of corresponding homogeneous spaces.
(Note that the rank of a homogeneous space $G/H$ is also
determined by the pair$(\g,\h)$, see~\cite{Pan1}.)

The method used in this paper is similar to the method of
I.~V.~Mykytyuk~\cite{mik1}. In particular, in~\cite{mik1} the
operation of coupling was called "the extension of a pair". Let us
remark that, using the operation of coupling of two pairs
"algebra-subalgebra", we do not need to consider the depth of a
subalgebra (M.~Brion), and may use induction only on the number
of simple components of the ambient algebra $\g$. The computation
of complexity of a homogeneous space is based on the study of
stationary subalgebras of general position (we use Elashvili's
Tables~\cite{El}) and Panyushev's formul\ae $\ $(Theorems 1 and
2).

The authors are grateful to D.~A.~Timashev for stimulative ideas
and numerous important remarks.


\section{Some facts about the complexity and the rank of homogeneous spaces}

Everywhere in the text Lie algebras of reductive groups are
denoted by the corresponding Gothic lowercase letters to the
exclusion of the exceptional simple Lie algebras. By tradition,
they are denoted by the Latin capital letters. In this section we
discuss some well-known facts about the complexity of homogeneous
spaces.

\begin{lemma}\label{lss}
Let $F$ be a reductive subgroup of $G$, and $H$ be an algebraic
subgroup of $F$. Then

1) \ $c_G(G/F)\le c_G(G/H)$;

2) \ $c_F(F/H)\le c_G(G/H)$.
\end{lemma}

\begin{proof}
By~\cite{vin1}, the complexity of a $G$-variety $X$
equals the modality of the induced action of $B$ on $X$. Let $X$
be the homogeneous space $G/H$. Then the complexity is the
modality of the double-sided action of $B\times H$ on $G$ or of
the right action of $H$ on $B\setminus G$. Hence, we have
inequality 1). In the same way we prove inequality 2). Indeed,
there exists a Borel subgroup $B$ of $G$ such that the
intersection of $B$ and $F$ is a Borel subgroup $B_F$ of $F$. Then
$B_F\setminus F$ is isomorphic to the $F$-orbit of the point $Be$
in $B\setminus G$.
\end{proof}

\begin{proposition}\label{os}
The following conditions are equivalent:

1) \ $c(G/H)=m$;

2) there exists an open subset $W\subset G$ such that for any
 $g\in W$ we have $\codim_{\g}(\h+\mathrm{Ad}(g)\bor)=m$.
\end{proposition}

\begin{proof}
Condition 1) holds if and only if  there exists an
open subset $V\subset G$ such that for any $g\in V$ the
codimension of $BgH$ (or of $g^{-1}BgH$) in $G$ is equal to $m$.
The latter is equivalent to condition 2) for $W=V^{-1}$.
\end{proof}

\medskip

We say that a pair $(\g, \h)$ has complexity $m$ if condition 2)
holds. In the sequel by $c(\g, \h)$ denote the complexity of a
pair $(\g, \h)$ . Proposition~1 allows to reduce the
classification of homogeneous spaces of complexity one to the
classification of corresponding pairs $(\g, \h)$.

\begin{proposition}
Let $\g=\g^s\oplus\cen$ be the decomposition into the direct sum
of the maximal semisimple ideal and the center. Then the
complexity of the pair $(\g, \h)$ equals the complexity of the
pair $(\g^s, pr_1(\h))$, where $pr_1: \g\to \g^s$ is the
projection along $\cen$.
\end{proposition}

\begin{proof}
Since $\cen$ is contained in $\mathrm{Ad}(g)\bor$
for any $g\in G$, we have
$\h+\mathrm{Ad}(g)\bor=pr_1(\h)+\mathrm{Ad}(g)\bor$.
\end{proof}

\medskip

Proposition~2 allows to consider only semisimple algebras $\g$. In the
sequel the subgroup $H$ and the subalgebra $\h$ are supposed to be
reductive.

Consider the restriction of the adjoint representation
$\mathrm{Ad}: G\to GL(\g)$ to the subgroup $H$. Since the
subalgebra $\h$ is an invariant subspace for this representation
and $H$ is reductive, there exists an $H$-invariant
subspace $\m$ of $\g$ such that $\g=\h\oplus \m$. The $H$-module
$\m$ can be naturally identified with the tangent space to the
affine homogeneous space $G/H$ at the point $eH$. The linear
representation $\rho: H\to GL(\m)$ is said to be {\it the isotropy
representation} of the pair $(G,H)$. Its differential is uniquely
defined by the pair $(\g, \h)$. By $\s$ denote a stationary
subalgebra of general position (s.s.g.p.) for the representation
$d\rho$. It is known that $\s$ is reductive ~\cite[Th.3]{Pan1}.
Let $\bor$ be a Borel subalgebra of $\g$, $\bor_{\s}$ be a Borel
subalgebra of $\s$, and $N(\g)$ be the number of positive roots in
the root system corresponding to the algebra  $\g$. The next
theorem shows  that the rank of $G/H$ is uniquely determined by
the pair $(\g,\h)$. This justifies the notation $\rk(\g,\h)$.

\begin{theorem}\cite[Th.3]{Pan1}
The following formul\ae $\,$ hold:
$$
c(G/H)= c(\g,\h)=\dim\,\bor -\rk\g-\dim\,\h+\dim\,\bor_{\s}=
N(\g)-\dim\,\h+N(\s)+\rk\s,
$$
$$
\rk(G/H)=\rk(\g,\h)=\rk\g-\rk\s.
$$
\end{theorem}

For our calculations we also need the following theorem.

\begin{theorem}\cite[Th.1.2]{Pan3}
Suppose that a subalgebra $\h$ is contained in a Levi subalgebra
$\lev$ of a parabolic subalgebra $\pp$ of the algebra $\g$, and
$S_1$ is a stabilizer of general position for the action of $H$ on
$\lev/\h$. Then

$$
  c(\g,\h)=c(\lev,\h)+c_{S_1}(\pp^u),
$$
where $\pp^u$ is the unipotent radical of the algebra $\pp$. In
particular, if $c(\lev,\h)=0$, then

$$
  c(\g,\h)=\dim\,\pp^u-\dim\,\bor_{\s_1}+\dim\,\bor_1,
$$
where $\bor_1$ is a s.s.g.p. for the action of $\bor_{\s_1}$ on
$\pp^u$.

\end{theorem}


\section{Reductive subalgebras of semisimple Lie algebras}

Let $\g=\g_1\oplus\dots\oplus\g_s$ be a decomposition of a
reductive Lie algebra $\g$ into the direct sum of the simple
ideals. We say that $\g_i$ are the components of the algebra $\g$.
If $\g$ is semisimple, then all $\g_i$ are simple
(noncommutative) Lie algebras. If $\g$ is reductive, then some
$\g_i$ can be one-dimensional.

Let $\h$ be a reductive subalgebra of $\g$, and $\h(i)$ be the
projection of $\h$ on $\g_i$ along the other components. Then
$\h(i)$ is a reductive subalgebra of $\g_i$. We say that the
subalgebra $\h$ is {\it completely decomposable} if
$\h=\h(1)\oplus\dots\oplus\h(s)$. In the general case, $\h$ is
contained in the completely decomposable subalgebra
$\h^r=\h(1)\oplus\dots\oplus\h(s)$, and $\h^r$ is the minimal
completely decomposable subalgebra containing $\h$.

In order to describe reductive subalgebras of $\g$ we have to do
two steps. The first is the description of completely decomposable
subalgebras (this is reduced to the description of reductive
subalgebras of simple algebras) and the second is the description
of reductive subalgebras such that the projection on any component
of the algebra is surjective.

\begin{lemma}\label{lll}
Let $\g=\g_1\oplus\dots\oplus\g_s$ be a semisimple Lie algebra,
and $\h$ be a reductive subalgebra of $\g$. Suppose that the
projection of $\h$ on every component $\g_i$ is surjective; then
$\h$ is semisimple, and any its component $\h_j$ embeds in
$\g_1\oplus\dots\oplus\g_s$ by the rule $x\to (pr_{j1}(x),\dots,
pr_{js}(x))$, where $pr_{ji}: \h_j \to \g_i$ is either an
isomorphism or the zero map (this means that $\h_j$ embeds
diagonally into the product of some mutually isomorphic $\g_i$).
For any $i$ there exists a unique non-zero map $pr_{ji}$.
\end{lemma}

\begin{proof}
Since the projection $pr_i:\h\to \g_i$ is
surjective, we see that ideals of $\h$ are taken to ideals of
$\g_i$. As the algebra $\g_i$ is simple, the image of any $\h_j$
is either $0$ or $\g_i$. This implies that there do not exist
commutative ideals of $\h$, i.e., $\h$ is semisimple.
Moreover, in the second case $pr_{ji}$ is an isomorphism. Finally,
the kernel of the projection $pr_i: \h\to\g_i$ is an ideal of
$\h$, therefore it is a sum of some components of the algebra
$\h$. Thus $pr_i$ determines an isomorphism between $\g_i$ and the
sum of some $\h_j$ such that $\h_j$ is not contained in the kernel
of $pr_i$. This proves the last statement of the Lemma.
\end{proof}

It is easy to generalize Lemma~\ref{lll} to a reductive $\g$. In
this case the center of $\h$ is the center of $\g$.

\smallskip

Let $\vv=\vv_1\oplus\dots\oplus\vv_t$ be a reductive Lie algebra
and $\phi_{i,j}:\vv_i\to \vv_j$ be an isomorphism for some $i\ne
j$. We say that the subalgebra
$$
\h=\{(x_1,\dots,x_i,\dots,x_{j-1},\phi_{ij}(x_i),x_{j+1},\dots,x_t)
\mid x_l\in \vv_l  \}
$$
is obtained from the algebra $\vv$ by {\it the operation of
coupling}. In general, $\h$ is obtained from $\vv$ by {\it a
multiple coupling} if it is obtained from $\vv$ by the successive
application of a finite number of couplings.

{\bf Conclusion.} Any reductive subalgebra $\h$ of a reductive Lie
algebra $\g=\g_1\oplus\dots\oplus\g_s$ can be obtained from a
completely decomposable subalgebra
$\ww=\ww(1)\oplus\dots\oplus\ww(s)$ ($=\h^r$), where $\ww(i)$ are
reductive subalgebras of $\g_i$, by a multiple coupling of some
components of the algebra $\ww$.


\section{The method of classification}
We say that a pair $(\g,\h)$ is {\it decomposable} if there exists
a non-trivial decomposition of $\g$ into the direct sum of ideals
$\g^1$ and $\g^2$ such that $\h=\h^1\oplus\h^2$, where
$\h^i=\g^i\cap\h$. By Proposition~\ref{os}, it follows that if
$\h=\h^1\oplus\h^2$ is a decomposable subalgebra of
$\g=\g^1\oplus\g^2$, then the complexity is additive:
$c(\g,\h)=c(\g^1,\h^1)+c(\g^2,\h^2)$. Thus, it is sufficient to
classify indecomposable pairs.

We shall make the classification by induction on the number $s$ of
components of the algebra $\g$. If $s=1$, then the algebra $\g$ is
simple and pairs of complexity~0 (respectively, 1) are represented
in Table 1 (respectively, in Table 2). Suppose that we have
classified all pairs $(\g, \h)$ of complexity $\le 1$ for all
algebras $\g$ such that $s=k$ (List $k$). Then we make a multiple
coupling in the subalgebra $\h^1\oplus\h^2$ of $\g^1\oplus\g^2$, where
$(\g^1,\h^1)$ is a pair from List $k$ and a pair $(\g^2,\h^2)$ is
from Table 1 or 2. Among obtained pairs we have all required pairs
$(\g, \h)$ with $s=k+1$. Moreover, one of the pairs $(\g^1,\h^1)$
or $(\g^2,\h^2)$ should be spherical (by Lemma~\ref{lss},~1), it
follows that $c(\g,\h)\ge c(\g^1, \h^1)+c(\g^2,\h^2)$). Note also
that it is sufficient to consider couplings of simple components
lying in different summands.

A pair $(\g,\h)$ is called {\it saturated} if the center of the
algebra $\h$ is a completely decomposable subalgebra of $\g$
(in~\cite{mik1} it was called "a main pair"). This means that $\h$
is obtained from a completely decomposable subalgebra of the
algebra $\g$ by couplings of noncommutative components. We
classify saturated indecomposable pairs of complexity one. As a
consequence, we obtain the classification of all indecomposable
pairs of complexity one, see Remark 2 below.

Let $(\g, \h)$ be a pair of complexity $\le 1$ and $\h_i$ be a
component of $\h$. We say that $\h_i$ is {\it valent}
(respectively, {\it 1-valent}), if $\h_i$ is not one-dimensional
and the subalgebra obtained by the coupling of the component
$\h_i$ in the subalgebra $\h\oplus\h_i$ is spherical in
$\g\oplus\h_i$ (respectively, of complexity 1). Such couplings are
called {\it elementary}.

The notion of valency allows to reduce the number of couplings
under consideration. Namely, suppose that we make a non-elementary
coupling. Then to obtain all desired pairs we should make
couplings of either two valent components or valent and 1-valent
components and the multiple couplings in the obtained algebras.
(Using Theorem~1, it is easy to show that the complexity of the
direct sum of two elementary couplings
$(\g_1\oplus\h\oplus\g_2\oplus\h, \h_1\oplus\h\oplus\h_2\oplus\h)$
is not greater than the complexity of the non-elementary coupling
$(\g_1\oplus\g_2, \h_1\oplus\h\oplus\h_2)$.)

\bigskip

{\it Step 1.}\ At first we determine all valent and 1-valent
components of pairs $(\g,\h)$, where $\g$ is simple. By the
definition of complexity, it follows that $c(\g,\h)\ge
N(\g)-\dim\h$. Consequently for the valency (respectively, the
1-valency) of the component $\h_i$ it is necessary that
$N(\g)+N(\h_i)-\dim\h\le 0$ (respectively, $\le 1$). By means of
this condition from pairs of Tables 1 and 2 we obtain the
following list of candidates for the valency or the 1-valency
(they are underlined):

$$ (\ssl_3, \underline{\so_3}),
(\ssl_{n+m}, \underline{\ssl_m}\oplus\ssl_n\oplus\cc)~\footnote{
under some conditions on $m,n$}, (\ssl_{n+m},
\underline{\ssl_m}\oplus\ssl_n)~\footnote{ under some conditions
on  $m,n$}, (\so_{n+m}, \underline{\so_m}\oplus\so_n)~\footnote{
under some conditions on  $m,n$},
$$
$$
(\so_7, \underline{G_2}),      
(\ssp_4, \underline{\ssl_2}\oplus\cc)~\footnote{ the embedding
corresponds to item 17 of Table 1}, (\ssp_{2(n+m)},
\underline{\ssp_{2m}}\oplus\ssp_{2n})~\footnote{ under some
conditions on $m,n$}, (G_2, \underline{\ssl_3}),  
$$
$$
(G_2, \underline{\ssl_2}\oplus\ssl_2),     
(F_4, \underline{\ssl_2}\oplus\ssp_6),     
(E_6, \underline{\ssl_2}\oplus\ssl_6),     
(E_7, \underline{\ssl_2}\oplus\so_{12}),   
$$
$$
(E_8, \underline{\ssl_2}\oplus E_7),         
(\ssp_{2n}, \underline{\ssl_2}\oplus\ssl_2\oplus\ssp_{2n-4}).
$$

To know which of these components are valent or 1-valent we use
Theorems 1 and 2.

\bigskip

{\it Remark 1.} \ Consider the pairs $(\g_1, \ssl_2\oplus\h_1)$
($\s_1$ is a s.s.g.p. for the isotropy representation, the
complexity equals $c_1$), $(\g_2, \ssl_2\oplus\h_2)$ ($\s_2$ is a
s.s.g.p., the complexity equals $c_2$) and the coupling of the
component $\ssl_2$ in the subalgebra: $(\g_1\oplus\g_2,
\ssl_2\oplus\h_1\oplus\h_2)$ ($\s$ is a s.s.g.p., the complexity
equals $c$). By Theorem 1, it follows that $c=c_1+c_2$ if and only
if in one of the pairs the projection of $\s_i$ on the component
$\ssl_2$ is surjective and for another component this projection
is non-zero. The condition $c=c_1+c_2+1$ holds if and only if
either one of the projections is surjective and another is zero or
the both projections are one-dimensional (their images are Cartan
subalgebras of $\ssl_2$). In other cases $c\ge c_1+c_2+2$. The
projection of the algebra $\s$ on the component $\ssl_2$ is zero
except for the case where the projections of $\s_1$ and $\s_2$ on
the component $\ssl_2$ are surjective. In this case the projection
of $\s$ on the component $\ssl_2$ is one-dimensional.

\medskip
By Remark~1, it follows that after an elementary coupling of
components isomorphic to $\ssl_2$ the complexity increases
by no more then 1 and
an $\ssl_2$-component of a spherical pair is always valent or
1-valent (in the first case the pair is contained in Table~3). By
the method of exception, we obtain 1-valent components from
spherical pairs:

$$
(\ssl_3, \underline{\so_3}), (G_2,
\underline{\ssl_2}\oplus\ssl_2), (F_4,
\underline{\ssl_2}\oplus\ssp_6),
$$
$$
(E_6, \underline{\ssl_2}\oplus\ssl_6), (E_7,
\underline{\ssl_2}\oplus\so_{12}), (E_8, \underline{\ssl_2}\oplus
E_7), (\ssp_4, \underline{\ssl_2}\oplus\cc),
$$
$$
(\ssl_3, \underline{\ssl_2}), (\so_{n+3},
\underline{\so_3}\oplus\so_n).
$$

It is easy to show that we have the next 1-valent components from
the candidates originated from the pairs of complexity one:

$$
(\ssp_{2n}, \underline{\ssl_2}\oplus\ssl_2\oplus\ssp_{2n-4}),
(\ssl_{4}, \underline{\ssl_2}\oplus\ssl_2).
$$

Let us pass to the question about the valency of components
different from $\ssl_2$. We apply Theorem 2 to the couplings
$(\ssl_{n+m}\oplus\ssl_m,
\underline{\ssl_m}\oplus\ssl_n\oplus\cc)$ and
$(\ssl_{n+m}\oplus\ssl_m, \underline{\ssl_m}\oplus\ssl_n)$. In the
first case $\bor_1=\ssc\oplus\bor^{n-m}$ (here $\bor^{n-m}$ is a
Borel subalgebra of $\ssl_{n-m}$) if $m<n$ and $\bor_1=0$ if $m\ge
n$. We obtain the valent components:
$$
(\ssl_{n+2}, \underline{\ssl_2}\oplus\ssl_n\oplus\cc),  
(\ssl_{m+1}, \underline{\ssl_m}\oplus\cc),             
(\ssl_{n+2}, \underline{\ssl_2}\oplus\ssl_n, \ n\geq 3),
$$
and the 1-valent components:
$$
(\ssl_{n+3}, \underline{\ssl_3}\oplus\ssl_n\oplus\cc, \ n\geq 2),
(\ssl_{n+3}, \underline{\ssl_3}\oplus\ssl_n, \ n\geq 4),
(\ssl_{m+1}, \underline{\ssl_m}), (\ssl_{4},
\underline{\ssl_2}\oplus\ssl_2).
$$

Consider the pair $(\so_{n+m}\oplus\so_m,
\underline{\so_m}\oplus\so_n)$. Here the isotropy representation
has the form $m=\p_1\otimes\p_1'\oplus ad(\so_m)$. Hence,
$\s=\so_{n-m}\subset\so_n$ if $n>m$ and $\s=0$ if $n\le m$. Thus
we obtain the valent component $(\so_{m+1},\underline{\so_m})$ and
the 1-valent component $(\so_{n+3},
\underline{\so_3}\oplus\so_n)$.

Consider a pair $(\ssp_{2(n+m)}\oplus\ssp_{2m},
\underline{\ssp_{2m}}\oplus\ssp_{2n})$. The isotropy
representation has the form $m=\p_1\otimes\p_1'\oplus
ad(\ssp_{2m})$. Thus we have $\s=\cc\oplus\ssp_{2n-2}$ if $m=1$
and $\s=\ssp_{2(n-m)}\subseteq\ssp_{2n}$ if $m\ge 2$ ($\s=0$ if
$n\le m$). Therefore we obtain the valent components
$(\ssp_{2(n+2)}, \underline{\ssp_{4}}\oplus\ssp_{2n})$,
$(\ssp_{2(n+1)}, \underline{\ssp_{2}}\oplus\ssp_{2n})$ and the
1-valent component $(\ssp_{8},
\underline{\ssp_{6}}\oplus\ssp_{2})$.

The isotropy representation of the coupling $(\so_7\oplus G_2,
\underline{G_2})$ has the form $m=\p_1\oplus ad(G_2)$. Here $\s=0$
and the component is 1-valent.

Consider the pair $(G_2\oplus\ssl_3, \underline{\ssl_3})$. Here
$\ssl_3$ embeds in $G_2$ "via the long roots". The short roots are
linearly independent on the Cartan subalgebra of $\ssl_3$, $\s=0$,
and the component is 1-valent.

Thus we obtain all valent and 1-valent simple components of pairs
$(\g, \h)$ of complexity $\le 1$, where $\g$ is simple.

\bigskip

{\it Step 2.}\ Now it is convenient to deviate from the standard
induction on the parameter $s$. We indicate all indecomposable
saturated pairs $(\g,\h)$ of complexity one such that $(\g,\h)$ is
obtained from pairs corresponding to simple algebras by couplings
of components isomorphic to $\ssl_2$.

By Remark 1 and the list of all valent and 1-valent components
isomorphic to  $\ssl_2$ (or by an explicit form of the isotropy
representation), it follows that if $(\g,\h)$ is a pair from the
list of candidates (Step 1) such that $\h$ has a component
isomorphic to $\ssl_2$, then the projection of the subalgebra $\s$
on this component is zero except for the pairs

1) $(\ssl_{n+2}, \ssl_2\oplus\ssl_n\oplus\cc), \ n\ge 1$, \
   $(\ssl_{n+2}, \ssl_2\oplus\ssl_n), \ n\ge 2$,\
   $(\ssp_{2n},  \ssl_2\oplus\ssl_2\oplus\ssp_{2n-4}), \ n\ge 3$~---
   here the projection is one-dimensional;

2) $(\ssp_{2(n+1)}, \ssl_2\oplus\ssp_{2n})$ ---
    here the projection is surjective.

Remark 1 and the information about the projections allow to find
all pairs of complexity 1 obtained by couplings of components
$\ssl_2$ from pairs $(\g,\h)$, where $\g$ is simple. By this way,
we obtain items 1 -- 22 of Table 4.

\bigskip

{\it Step 3.}\ By direct calculations based on examination of the
dimension and the isotropy representation, it is easy to prove
that all non-elementary couplings of valent and 1-valent
components different from $\ssl_2$ lead to pairs of
complexity~$>1$

Finally we put items 24 -- 29 corresponding to obtained 1-valent
components to Table 4-4.

\bigskip

{\it Step 4.}\ Now we enumerate valent and 1-valent components
different from $\ssl_2$ for pairs of complexity $\le 1$ with
$s=2$. Consider a pair $(\g,\h)$. Suppose that we obtain this pair
by couplings in the pairs $(\g_i,\h_i)$, where $\g_i$ are simple;
then we have to consider only components corresponding to the
valent and the 1-valent components of $(\g_i,\h_i)$.

Start with item 9 of Table 3. The component $\h$ of this pair is
1-valent if and only if $N(\h)-\rk\h=1$. This holds only if
$\h=\ssl_3$ (this corresponds to item 30 of Table 4-4). By Lemma
~\ref{lss},~2), it follows that if $\h_i$ is a component of a pair
of complexity one, $\h_i\ne \ssl_2, \ssl_3$ then $\h_i$ can not be
embedded in more than two simple components of the algebra $\g$.

It remains to check the valency or the 1-valency of the underlined
components  of the next pairs:

$$
(\g_1\oplus\ssp_{2n+4},
\h_1\oplus\ssp_{2n}\oplus\underline{\ssp_4}), \ (n=1,2),
$$
$$
(\g_1\oplus\ssp_8, \h_1\oplus\ssl_2\oplus\underline{\ssp_6}),
$$
$$
(\g_1\oplus\ssl_{n+3}, \h_1\oplus\ssl_n\oplus\underline{\ssl_3}),
\ (n=2,3),
$$
$$
(\g_1\oplus\ssl_{n+3}, \h_1\oplus\gl_n\oplus\underline{\ssl_3}), \
(n=2,3),
$$
$$
(\g_1\oplus\g_2, \h_1\oplus\underline{\ssl_3}\oplus\h_2),
$$
where the first component of $\h$ is embedded only in $\g_1$, the
third one only in the second component of $\g$, and the  second
one in the both components of $\g$.

Theorem 1 shows that a necessary condition for underlined
components to be valent or 1-valent is that the projection of the
algebra $\s$ on this component for the indicated pair contains a
Borel subalgebra of dimension no less than 3, 8, 2, 2 ¨ 2
respectively. Using the explicit form of a subalgebra $\s$
computed for all pairs of Table 4 for $s=2$ we have that this
condition holds only for $(\ssp_{2n+2}\oplus\ssp_6,
\ssp_{2n}\oplus\ssl_2\oplus\ssp_4)$. It is easy to show that the
component $\ssp_4$ is 1-valent here (Table~4-4, item 23).

\bigskip

{\it Step 5.}\ For $s=3$ we have to consider a non-elementary
coupling of the components different from $\ssl_2$ only for
item~16 of Table~1 with $m=2$ and item 3 of Table 3 with $m=2$
(here the joint component is $\ssp_4$). The complexity of this
coupling is greater than 3.

\bigskip

{\it Step 6.}\ Similarly to Step~4 it can be checked that a pair
with $s=3$ has no valent and 1-valent components different from
$\ssl_2$.

\bigskip

It remains to note that for any pair from Table 4 if there exist
two different isomorphisms determining embeddings of a component
of the subalgebra into a component of the algebra, then these
isomorphisms give the subalgebras which are taken to each other by
an automorphism of the ambient algebra. Thus we proved

\begin{theorem}
 All indecomposable saturated pairs $(\g,\h)$ of complexity
one, where $\g$ is a simple Lie algebra and $\h$ is its reductive
subalgebra, are represented in  Table 4 (up to automorphisms of
the algebra $\g$).

\end{theorem}

It is interesting to note that in Table 4 there are no pairs of
even rank greater than 8 and there is only one series of pairs of
odd rank greater than 7 (item 27).

\bigskip

{\it Remark 2.} \ We shall show how to describe all indecomposable
pairs of complexity one having the description of all saturated
indecomposable pairs of complexity one. The idea of this method
was suggested by D.~A.~Timashev. For any indecomposable pair
$(\g,\h)$ we have the saturated (maybe decomposable) pair
$(\g,\tilde\h)$, where the semisimple components of $\h$ and
$\tilde\h$ are the same (we denote them by $\h^s$) and the center
of $\tilde\h$ (we denote it by
$\tilde\z=\z_1\oplus\dots\oplus\z_s$) is the minimal completely
decomposable subalgebra containing the center of $\h$ (we denote
the center of $\h$ by $\z$). If $c(\g,\h)=1$, then
$c(\g,\tilde\h)\le 1$. In the algebra $\tilde\z$ there exists the
largest subalgebra $\ppp$ such that $$
c(\g,\h^s)=c(\g,\h^s+\ppp)=c(\g,\tilde\h)-(\dim\tilde\z-\dim\ppp).
$$ In fact, consider the maximal subtorus of the center of the group $\tilde H$
such that the generic $H^s$-orbit on $G/B$ is invariant under its
action. The subalgebra $\ppp$ is the tangent algebra to this
subtorus. Note that $\ppp$ can be decomposed into the direct sum
of subalgebras corresponding to the decomposition of the pair
$(\g,\tilde\h)$ into indecomposable components. Thus it is
sufficient to find the subalgebra $\ppp$ for every saturated
indecomposable pair. Tables 1, 2, 3, 4 imply that the subalgebra
$\ppp$ is zero in all
 cases except for the followings:

1)\ Table 1: \ it.2, $n\ne m$;
              \ it.6;
              \ it.7, $n$ is odd;
              \ it.25;

    Table 3: \ it.1, $n\ge 3$;

   Table 4: \ it.8, $n\ge 3,\ m\ge 3$;
              \ it.9, $n\ge 3$;
              \ it.16, $n\ge 3$;
              \ it.25, $n\ge 4$
--- in all these cases $\ppp=\tilde\z$;

2)\ Table 2, it.2: here $\ppp$ is a one-dimensional subalgebra
described in it.3;

3)\ Table 4, it.8, $n=1,2, m\ge 3$ and $m=1,2, n\ge 3$: here
$\ppp$ is one-dimensional and corresponds to the bigger component.

It remains to note that the pair $(\g,\h)$ has the complexity 1 if
and only if either $c(\g,\tilde\h)=1$ and the restriction of the
projection $\tilde\z \to \tilde\z/\ppp$ to the subalgebra $\z$ is
 surjective, or $c(\g,\tilde\h)=0$  and the image of $\z$ under the projection
is a hyperplane.

\medskip

{\bf Notations and conventions used in the Tables}. The indexes
$n,m$ in Tables 1 and 2 are assumed to be $\ge 1$ and the indexes
$n,m,k,l$ in Tables 3 and 4 are assumed to be $\ge 0$ (if the
converse is not mentioned). If the index of a classical algebra
has a non-positive value (for example, $\ssp_{2n-2}$ with
$n=0,1$), then the corresponding algebra is assumed to be zero.
The symbol $\delta_i^j$ equals 1 if $i=j$ and 0 otherwise. In
Tables 1 and 2 the column "$\h\subset\g$" determines the embedding
of the subalgebra $\h$ in the algebra $\g$ by the restriction of
the first fundamental representation of $\g$ to $\h$ (the
numeration of the fundamental weights is taken from~\cite{vo}).
By $\pi_i$ and $\pi_i'$ we denote the fundamental weights of the
first and the second simple components of the algebra $\h$
respectively. By $\epsilon$ we denote the fundamental weight of
the one-dimensional central subalgebra $\cc$. We use
multiplicative notation for weights (for example, $\pi_1^2$ is the
highest weight with the first mark 2 and other marks zero). By the
sing "+" we denote the direct sum of representations. By the
symbol 1 we denote the one-dimensional trivial representation, by
2 --- the two-dimensional trivial representation, and so on. In
Tables 3 and 4 segments denote the embedding of the subalgebra in
the algebra. If a segment does not determine the embedding
uniquely, then the embedding is indicated on the segment as in
Tables 1 and 2.

\pagebreak

\centerline{\large Table 1}

\vspace{0.7cm}

\footnotesize

\centerline{ \small
\begin{tabular}{|c|c|c|c|c|}
 \hline
  & $\g$ & $\h$ & $\h\subset \g$ & $\m$\\
\hline
 1 & $\ssl_n, \ n\ge 2 $ & $\so_n$  & $\p_1$ & $\pi_1^2$ \ $(n\ne 2,4)$ \\
\hline
 2 & $\ssl_{n+m}$ & $\ssl_n\oplus\ssl_m\oplus\ssc$ &
$\p_1\e^m+\p'_1\e^{-n}$
 & $\p_1\otimes\p_{m-1}'\e^{m+n} +\p_{n-1}\otimes\p_1'\e^{-m-n}$\\
\hline
 3 & $\ssl_{n+m}$, & $\ssl_n\oplus \ssl_m$ & $\p_1+\p'_1$
&$\p_1 \otimes \p_{m-1}' +\p_{n-1} \otimes \p_1'+ 1$ \\
  & $m\ne n$ & & & \\
\hline
 4 & $\ssl_{2n},  \ n\ge 2$ & $\ssp_{2n}$  & $\p_1$ & $\p_2$ \\
\hline
 5 & $\ssl_{2n+1}$  & $\ssp_{2n}$ & $\p_1+1$&  $\p_1^2+\p_2+1$ \\
\hline
 6 & $\ssl_{2n+1}$ & $\ssp_{2n}\oplus\ssc$  & $\p_1\e+\e^{-2n}$ &
$\p_1\e^{2n+1}+\p_1\e^{-2n-1}+\p_2$ \\
\hline
 7 &$\so_{2n}$ & $\ssl_n\oplus\ssc$ &
$\p_1\e+\p_{n-1}\e^{-1}$&$\p_2\e^2+\p_{n-2}\e^{-2}$\\
\hline
 8 & $\so_{4n+2}$  & $\ssl_{2n+1}$ & $\p_1+\p_{2n}$ & $\p_2+\p_{2n-1}+1$\\
\hline
 9 & $\so_{2n+1}$ & $\ssl_n\oplus\ssc$  & $\p_1\e+\p_{n-1}\e^{-1}+1$ &
$\p_1\e+\p_2\e^2+\p_{n-2}\e^{-2}+\p_{n-1}\e^{-1}$ \\
\hline 10 & $\so_{n+m}$ & $\so_n\oplus \so_m$ \ & $\p_1+\p_1'$&
$\p_1\otimes\p_1'$
\\
\hline
11 & $\so_9$ & $\so_7$  & $\p_3+1$ &$\p_1+\p_3$\\
\hline
12 & $\so_7$ & $G_2$  & $\p_1$ &$\p_1$\\
\hline
13 & $\so_8$ & $G_2$  & $\p_1+1$& $\p_1^2$\\
\hline 14 & $\so_{10}$ & $\so_7\oplus\ssc$  & $\p_3+ \e+\e^{-1}$ &
$\pi_1+\pi_3\e+\pi_3\e^{-1}$ \\
\hline 15 & $\ssp_{2n}$ & $\ssl_n\oplus\ssc$  &
$\p_1\e+\p_{n-1}\e^{-1}$ &
$\pi_1^2\e^2+\pi_{n-1}^2\e^{-2}$ \\
\hline 16 & $\ssp_{2(n+m)}$  & $\ssp_{2n}\oplus \ssp_{2m}$ &
$\p_1+\p_1'$&
 $\p_1\otimes\p_1'$\\
\hline 17 & $\ssp_{2n}$ & $\ssp_{2(n-1)}\oplus\ssc$ &
$\p_1+\e+\e^{-1}$&
$\p_1\e+\p_1\e^{-1}+\e^2+\e^{-2}$ \\
\hline
18 & $G_2$ & $\ssl_3$  & $\p_1+\p_2+1$& $\p_1+\p_2$ \\
\hline 19 & $G_2$  & $\ssl_2\oplus \ssl_2$ &
$\p_1^2+\p_1\otimes\p_1'$ &
$\p_1^3\otimes\p_1'$\\
\hline
20 & $F_4$ & $\so_9$  & $\p_1+\p_4+1$ &$\p_4$\\
\hline 21 & $F_4$ & $\ssp_6\oplus \ssl_2$  &
$\p_2+\p_1\otimes\p_1'$ &
$\p_3\otimes\p_1'$ \\
\hline
22 & $E_6$ & $\ssp_8$  & $\p_2$ & $\p_4$ \\
\hline
23 & $E_6$ & $F_4$  & $\p_1+1$ & $\p_1$\\
\hline
24 & $E_6$  & $\so_{10}$ & $\p_1+\p_5+1$ &$\p_4+\p_5+1$\\
\hline 25 & $E_6$& $\so_{10}\oplus\ssc$  &
$\p_1\e^2+\p_5\e^{-1}+\e^{-4}$ &
$\p_4\e^6+\p_5\e^{-6}$ \\
\hline 26 & $E_6$  & $\ssl_6\oplus \ssl_2$ &
$\p_4+\p_1\otimes\p_1'$ &
$\p_3\otimes\p_1'$\\
\hline 27 & $E_7$  & $E_6\oplus\ssc$ &
$\p_1\e+\p_5\e^{-1}+\e^{3}+\e^{-3}$&
$\p_1\e^4+\p_5\e^{-4}$ \\
\hline
28 & $E_7$  & $\ssl_8$ & $\p_2+\p_6$ & $\p_4$ \\
\hline 29 & $E_7$ & $\so_{12}\times \ssl_2$ &
$\p_6+\p_1\otimes\p_1'$ &
$\p_5\otimes\p_1'$\\
\hline
30 & $E_8$ & $\so_{16}$  & $\p_2+\p_8$ & $\p_7$ \\
\hline 31 & $E_8$& $\ssl_2\times E_7$ &
$\p_1\otimes\p_1'+\p_1^2+\p_6'$ &
$\p_1\otimes\p_1'$\\
\hline
\end{tabular}
}

\normalsize \vspace{1cm} \

\footnotesize

\pagebreak

\centerline{\large Table 2}

\vspace{0.7cm}

\centerline{
\begin{tabular}{|c|c|c|c|c|}
\hline
  & $\g$ & $\h$ & $\h\subset \g$ & $\m$\\
\hline
 1& $\ssl_{2n}$ & $\ssl_n\oplus \ssl_n$  & $\p_1+\p_1'$&
$\p_1\otimes \p_1'+\p_{n-1}\otimes\p_{n-1}'+1$ \\
\hline
 2 &$\ssl_n$& $\ssl_{n-2}\oplus \ssc \oplus \ssc$ &
$\p_1\e_1\e_2+\e_1^{2-n}+\e_2^{2-n}$&
$\p_1\e^{n-1}_1\e_2+\p_1\e_1\e_2^{n-1}+$\\
& $n\geq 3$&&&$+\e_1^{2-n}\e_2^{n-2}+\e_1^{n-2}\e_2^{2-n}+$\\
&&&&$+\p_{n-3}\e^{1-n}_1\e_2^{-1}+\p_{n-3}\e_1^{-1}\e_2^{1-n}$\\
\hline 3 & $\ssl_n$& $\ssl_{n-2}\oplus \ssc$  &
$\p_1\e+\e^{d_1}+\e^{d_2}$ &
$\p_1\e^{1-d_1}+\p_1\e^{1-d_2}+$\\
   & $n\ge 5$&  & $d_1 \ne d_2$, $d_1+d_2=2-n$&
$+\e^{d_1-d_2}+\e^{d_2-d_1}+$\\
&&&&$+\p_{n-3}\e^{d_1-1}+\p_{n-3}\e^{d_2-1}$\\
\hline
 4 & $\ssl_6$& $\ssp_4\oplus \ssl_2 \oplus \ssc$  & $\p_1\e+\p_1'\e^{-2}$
&
$\p_1\otimes\p_1'\e+\p_2+$\\
&&&&$+\p_1\otimes\p_1'\e^{-1}$\\
\hline
 5& $\so_n$ & $\so_{n-2}$  & $\p_1+2$ &$\p_1+\p_1+1$\\
   &     $n\ge 5$&& &\\
\hline
 6& $\so_{2n+1}$ & $\ssl_n$  & $\p_1+\p_{n-1}+1$ &$\p_1+\p_2+\p_{n-1}+$\\
   &  $n\ge 3$ && & + $\p_{n-2}+1$\\
\hline
 7& $\so_{4n}$ & $\ssl_{2n}$  & $\p_1+\p_{2n-1}$ &$\p_2+\p_{2n-2}$\\
   &     $n\ge 2$ & & &\\
\hline
 8& $\so_9$  & $G_2\oplus \ssc$ & $\p_1+\e+\e^{-1}$
&$\p_1+\p_1\e+\p_1\e^{-1}$\\
\hline
 9& $\so_{11}$  & $\ssl_2\oplus \so_7$ & $\p_1^2+\p_3'$
& $\p_1^2\otimes\p_3'+\p_1'$\\
\hline
 10& $\so_{10}$ & $\so_7$  & $\p_3+2$ &$\p_1+\p_3+\p_3+1$\\
\hline
 11& $\ssp_{2n}$ & $\ssl_n$  & $\p_1+\p_{n-1}$ &$\p_1^2+\p_{n-1}^2+1$\\
\hline
 12& $\ssp_{2n}$  & $\ssp_{2n-2}$ & $\p_1+2$ &$\p_1+\p_1+3$\\
              & $n\ge 2$ && &\\
\hline
 13 & $\ssp_{2n}$& $\ssp_{2n-4}\oplus \ssl_2\oplus \ssl_2$  &
$\p_1+\p_1'+\p_1''$& $\p_1\otimes\p_1'+\p_1\otimes\p_1''+$\\
&                                   $n\ge 3$ && &$+\p_1'\otimes\p_1''$\\
\hline
14& $\ssp_4$ &  $\ssl_2$  & $\p_1^3$ &$\p_1^6$\\
\hline 15 & $E_6$& $\so_9\oplus \ssc$  &
$\p_1\e^2+\e^2+\p_4\e^{-1}+\e^{-4}$ &
$\p_1+\p_4\e+\p_4\e^{-1}$\\
\hline
16  &$E_7$  & $E_6$ &  $\p_1+\p_5+2$ &$\p_1+\p_5+1$\\
\hline
17& $F_4$ & $\so_8$  & $\p_1+\p_3+\p_4+2$ &$\p_1+\p_3+\p_4$\\
\hline
\end{tabular}
}

\normalsize \vspace{1cm}

\pagebreak

\centerline{ \large Table 3}

\vspace{0.7cm} \centerline{ \footnotesize
\begin{tabular}{|c|c|c|c|}
\hline
 & $\h\subset \g$ & $\s$ & rk\\
\hline 1  & \begin{picture}(120,60)
     \put(20,51){$\ssl_{n+2}$}
     \put(60,51){$\ssp_{2m+2}$}
     \put(10,15){\line(2,3){20}}
     \put(50,15){\line(-2,3){20}}
     \put(50,15){\line(2,3){20}}
     \put(90,15){\line(-2,3){20}}
     \put(10,15){\circle*{2}}
     \put(50,15){\circle*{2}}
     \put(90,15){\circle*{2}}
     \put(30,45){\circle*{2}}
     \put(70,45){\circle*{2}}
     \put(0,3){$\gl_{n}$}
     \put(40,3){$\ssl_2$}
     \put(80,3){$\ssp_{2m}$}
     \end{picture} $n\ge 1$ & $\gl_{n-2}\oplus\ssp_{2m-2}$ &
 $5-\delta_m^0-\delta_n^1$   \\
\hline 2  & \begin{picture}(120,60)
     \put(20,51){$\ssl_{n+2}$}
     \put(60,51){$\ssp_{2m+2}$}
     \put(10,15){\line(2,3){20}}
     \put(50,15){\line(-2,3){20}}
     \put(50,15){\line(2,3){20}}
     \put(90,15){\line(-2,3){20}}
     \put(10,15){\circle*{2}}
     \put(50,15){\circle*{2}}
     \put(90,15){\circle*{2}}
     \put(30,45){\circle*{2}}
     \put(70,45){\circle*{2}}
     \put(0,3){$\ssl_{n}$}
     \put(40,6){$\ssl_2$}
     \put(80,6){$\ssp_{2m}$}
     \end{picture}  $n\ge 3$ & $\ssl_{n-2}\oplus\ssp_{2m-2}$
& $6-\delta_m^0$ \\
\hline 3  & \begin{picture}(120,60)
     \put(20,51){$\ssp_{2n+2}$}
     \put(60,51){$\ssp_{2m+2}$}
     \put(10,15){\line(2,3){20}}
     \put(50,15){\line(-2,3){20}}
     \put(50,15){\line(2,3){20}}
     \put(90,15){\line(-2,3){20}}
     \put(10,15){\circle*{2}}
     \put(50,15){\circle*{2}}
     \put(90,15){\circle*{2}}
     \put(30,45){\circle*{2}}
     \put(70,45){\circle*{2}}
     \put(0,3){$\ssp_{2n}$}
     \put(40,3){$\ssl_2$}
     \put(80,3){$\ssp_{2m}$}
     \end{picture} &$\ssp_{2n-2}\oplus\ssp_{2m-2}\oplus\ssc$
& $3-\delta_m^0-\delta_n^0$  \\
\hline 4  & \begin{picture}(120,60)
     \put(0,51){$\ssp_{2n+2}$}
     \put(60,51){$\ssp_{2m+2}$}
     \put(90,51){$\ssp_{2l+2}$}
     \put(7,15){\line(0,1){30}}
     \put(67,15){\line(0,1){30}}
     \put(97,15){\line(0,1){30}}
     \put(7,15){\circle*{2}}
     \put(7,45){\circle*{2}}
     \put(67,15){\circle*{2}}
     \put(67,45){\circle*{2}}
     \put(97,15){\circle*{2}}
     \put(97,45){\circle*{2}}
     \put(0,3){$\ssp_{2n}$}
     \put(60,3){$\ssp_{2m}$}
     \put(90,3){$\ssp_{2l}$}
     \put(37,15){\line(-1,1){30}}
     \put(37,15){\line(1,1){30}}
     \put(37,15){\line(2,1){60}}
     \put(37,15){\circle*{2}}
     \put(30,3){$\ssl_2$}
     \end{picture} &$\ssp_{2n-2}\oplus\ssp_{2m-2}\oplus\ssp_{2l-2}$&
$6-\delta_m^0-\delta_n^0-\delta_l^0$ \\
\hline 5  & \begin{picture}(120,60)
     \put(15,51){$\ssp_{2n+2}$}
     \put(45,51){$\ssp_4$}
     \put(75,51){$\ssp_{2m+2}$}
     \put(7,15){\line(1,2){15}}
     \put(37,15){\line(-1,2){15}}
     \put(37,15){\line(1,2){15}}
     \put(67,15){\line(-1,2){15}}
     \put(67,15){\line(1,2){15}}
     \put(97,15){\line(-1,2){15}}
     \put(7,15){\circle*{2}}
     \put(37,15){\circle*{2}}
     \put(67,15){\circle*{2}}
     \put(22,45){\circle*{2}}
     \put(52,45){\circle*{2}}
     \put(97,15){\circle*{2}}
     \put(82,45){\circle*{2}}
     \put(0,3){$\ssp_{2n}$}
     \put(30,3){$\ssl_2$}
     \put(60,3){$\ssl_2$}
     \put(90,3){$\ssp_{2m}$}
     \end{picture}&$\ssp_{2n-2}\oplus\ssp_{2m-2}$&
$6-\delta_m^0-\delta_n^0$ \\
\hline 6   & \begin{picture}(80,60)
    \put(0,51){$\ssp_{2n+4}$}
    \put(40,51){$\ssp_4$}
    \put(10,15){\line(0,1){30}}
    \put(50,15){\line(-4,3){40}}
    \put(50,15){\line(0,1){30}}
    \put(10,15){\circle*{2}}
    \put(10,45){\circle*{2}}
    \put(50,15){\circle*{2}}
    \put(50,45){\circle*{2}}
    \put(0,3){$\ssp_{2n}$}
    \put(40,3){$\ssp_4$}
    \end{picture} $n\ge 1$ &$\ssp_{2n-4}$ &
$6-\delta_n^1$ \\
\hline 7   & \begin{picture}(80,60)
    \put(0,51){$\ssl_n$}
    \put(40,51){$\ssl_{n+1}$}
    \put(10,15){\line(0,1){30}}
    \put(10,15){\line(4,3){40}}
    \put(50,15){\line(0,1){30}}
    \put(10,15){\circle*{2}}
    \put(10,45){\circle*{2}}
    \put(50,15){\circle*{2}}
    \put(50,45){\circle*{2}}
    \put(0,3){$\ssl_n$}
    \put(40,3){$\ssc$}
    \end{picture} $n\ge 2$ & 0 & $2n-1$\\
\hline 8  & \begin{picture}(80,60)
     \put(0,51){$\so_n$}
     \put(40,51){$\so_{n+1}$}
     \put(30,15){\line(-2,3){20}}
     \put(30,15){\line(2,3){20}}
     \put(30,15){\circle*{2}}
     \put(10,45){\circle*{2}}
     \put(50,45){\circle*{2}}
     \put(20,3){$\so_n$}
     \end{picture} $n\ge 3$ & 0& $n$\\
\hline 9  & $\h$ is simple
     \begin{picture}(80,60)
     \put(10,51){$\h$}
     \put(50,51){$\h$}
     \put(30,15){\line(-2,3){20}}
     \put(30,15){\line(2,3){20}}
     \put(30,15){\circle*{2}}
     \put(10,45){\circle*{2}}
     \put(50,45){\circle*{2}}
     \put(30,3){$\h$}
     \end{picture} & Cartan subalgebra & rk $\h$\\
\hline
\end{tabular}
} \pagebreak

\centerline{\large Table 4-1}

\vspace{0.7cm}

\centerline{ \footnotesize
\begin{tabular}{|c|c|c|c|}
\hline
 & $\h\subset \g$ &$\s$& rk\\ \hline 1   &
\begin{picture}(80,66)
    \put(0,51){$\ssp_{2n+2}$}
    \put(40,51){$\ssl_3$}
    \put(10,15){\line(0,1){30}}
    \put(50,15){\line(-4,3){40}}
    \put(50,15){\line(0,1){30}}
    \put(10,15){\circle*{2}}
    \put(10,45){\circle*{2}}
    \put(50,15){\circle*{2}}
    \put(50,45){\circle*{2}}
    \put(54,30){$\p_1^2$}
    \put(0,3){$\ssp_{2n}$}
    \put(40,3){$\ssl_2$}
    \end{picture}&$\ssp_{2n-2}$& $4-\delta_n^0$ \\
\hline 2  & \begin{picture}(100,66)
     \put(20,51){$\ssp_{2n+2}$}
     \put(64,51){$G_2$}
     \put(10,15){\line(2,3){20}}
     \put(50,15){\line(-2,3){20}}
     \put(50,15){\line(2,3){20}}
     \put(90,15){\line(-2,3){20}}
     \put(10,15){\circle*{2}}
     \put(50,15){\circle*{2}}
     \put(90,15){\circle*{2}}
     \put(30,45){\circle*{2}}
     \put(70,45){\circle*{2}}
     \put(0,3){$\ssp_{2n}$}
     \put(40,3){$\ssl_2$}
     \put(84,3){$\ssl_2$}
     \end{picture} &$\ssp_{2n-2}$& $4-\delta_n^0$ \\
\hline 3  & \begin{picture}(100,66)
     \put(20,51){$\ssp_{2n+2}$}
     \put(64,51){$F_4$}
     \put(10,15){\line(2,3){20}}
     \put(50,15){\line(-2,3){20}}
     \put(50,15){\line(2,3){20}}
     \put(90,15){\line(-2,3){20}}
     \put(10,15){\circle*{2}}
     \put(50,15){\circle*{2}}
     \put(90,15){\circle*{2}}
     \put(30,45){\circle*{2}}
     \put(70,45){\circle*{2}}
     \put(0,3){$\ssp_{2n}$}
     \put(40,3){$\ssl_2$}
     \put(84,3){$\ssp_6$}
     \end{picture} & $\ssp_{2n-2}$ & $6-\delta_n^0$ \\
\hline 4  & \begin{picture}(100,66)
     \put(20,51){$\ssp_{2n+2}$}
     \put(64,51){$E_6$}
     \put(10,15){\line(2,3){20}}
     \put(50,15){\line(-2,3){20}}
     \put(50,15){\line(2,3){20}}
     \put(90,15){\line(-2,3){20}}
     \put(10,15){\circle*{2}}
     \put(50,15){\circle*{2}}
     \put(90,15){\circle*{2}}
     \put(30,45){\circle*{2}}
     \put(70,45){\circle*{2}}
     \put(0,3){$\ssp_{2n}$}
     \put(40,3){$\ssl_2$}
     \put(84,3){$\ssl_6$}
     \end{picture} & $\ssp_{2n-2}\oplus\cc\oplus\cc$ & $6-\delta_n^0$ \\
\hline 5  & \begin{picture}(100,66)
     \put(20,51){$\ssp_{2n+2}$}
     \put(64,51){$E_7$}
     \put(10,15){\line(2,3){20}}
     \put(50,15){\line(-2,3){20}}
     \put(50,15){\line(2,3){20}}
     \put(90,15){\line(-2,3){20}}
     \put(10,15){\circle*{2}}
     \put(50,15){\circle*{2}}
     \put(90,15){\circle*{2}}
     \put(30,45){\circle*{2}}
     \put(70,45){\circle*{2}}
     \put(0,3){$\ssp_{2n}$}
     \put(40,3){$\ssl_2$}
     \put(84,3){$\so_{12}$}
     \end{picture} & $\ssp_{2n-2}\oplus\ssl_2\oplus\ssl_2\oplus\ssl_2$
& $6-\delta_n^0$ \\
\hline 6  & \begin{picture}(100,66)
     \put(20,51){$\ssp_{2n+2}$}
     \put(64,51){$E_8$}
     \put(10,15){\line(2,3){20}}
     \put(50,15){\line(-2,3){20}}
     \put(50,15){\line(2,3){20}}
     \put(90,15){\line(-2,3){20}}
     \put(10,15){\circle*{2}}
     \put(50,15){\circle*{2}}
     \put(90,15){\circle*{2}}
     \put(30,45){\circle*{2}}
     \put(70,45){\circle*{2}}
     \put(0,3){$\ssp_{2n}$}
     \put(40,3){$\ssl_2$}
     \put(84,3){$E_7$}
     \end{picture} & $\ssp_{2n-2}\oplus\so_8$ & $6-\delta_n^0$ \\
\hline 7  & \begin{picture}(100,66)
     \put(20,51){$\ssp_{2n+2}$}
     \put(64,51){$\so_{k+3}$}
     \put(10,15){\line(2,3){20}}
     \put(50,15){\line(-2,3){20}}
     \put(50,15){\line(2,3){20}}
     \put(90,15){\line(-2,3){20}}
     \put(10,15){\circle*{2}}
     \put(50,15){\circle*{2}}
     \put(90,15){\circle*{2}}
     \put(30,45){\circle*{2}}
     \put(70,45){\circle*{2}}
     \put(0,3){$\ssp_{2n}$}
     \put(40,3){$\ssl_2=\so_3$}
     \put(84,3){$\so_k$}
     \end{picture} $k\ge 2$ & $\ssp_{2n-2}\oplus\so_{k-3}$ &
$5-\delta_n^0-\delta_k^2$\\
\hline 8  & \begin{picture}(100,66)
     \put(20,51){$\ssl_{n+2}$}
     \put(64,51){$\ssl_{m+2}$}
     \put(10,15){\line(2,3){20}}
     \put(50,15){\line(-2,3){20}}
     \put(50,15){\line(2,3){20}}
     \put(90,15){\line(-2,3){20}}
     \put(10,15){\circle*{2}}
     \put(50,15){\circle*{2}}
     \put(90,15){\circle*{2}}
     \put(30,45){\circle*{2}}
     \put(70,45){\circle*{2}}
     \put(0,3){$\gl_{n}$}
     \put(40,3){$\ssl_2$}
     \put(84,3){$\gl_m$}
     \end{picture} $n,m\ge 1$& $\gl_{n-2}\oplus\gl_{m-2}$ &
$6-\delta_n^1-\delta_m^1$ \\
\hline
\end{tabular}
} \pagebreak

\centerline{\large Table 4-2}

\vspace{0.7cm}

\centerline{
\begin{tabular}{|c|c|c|c|}
\hline
 & $\h\subset \g$ & $\s$ &rk\\
\hline 9  & \begin{picture}(100,66)
     \put(20,51){$\ssl_{n+2}$}
     \put(64,51){$\ssl_{m+2}$}
     \put(10,15){\line(2,3){20}}
     \put(50,15){\line(-2,3){20}}
     \put(50,15){\line(2,3){20}}
     \put(90,15){\line(-2,3){20}}
     \put(10,15){\circle*{2}}
     \put(50,15){\circle*{2}}
     \put(90,15){\circle*{2}}
     \put(30,45){\circle*{2}}
     \put(70,45){\circle*{2}}
     \put(0,3){$\gl_{n}$}
     \put(40,3){$\ssl_2$}
     \put(84,3){$\ssl_m$}
     \end{picture} $m\ge 3$,$n\ge 1$ & $\gl_{n-2}\oplus\ssl_{m-2}$
& $7-\delta_n^1$ \\
\hline 10  & \begin{picture}(100,66)
     \put(20,51){$\ssl_{n+2}$}
     \put(64,51){$\ssl_{m+2}$}
     \put(10,15){\line(2,3){20}}
     \put(50,15){\line(-2,3){20}}
     \put(50,15){\line(2,3){20}}
     \put(90,15){\line(-2,3){20}}
     \put(10,15){\circle*{2}}
     \put(50,15){\circle*{2}}
     \put(90,15){\circle*{2}}
     \put(30,45){\circle*{2}}
     \put(70,45){\circle*{2}}
     \put(0,3){$\ssl_{n}$}
     \put(40,3){$\ssl_2$}
     \put(82,3){$\ssl_m$}
     \end{picture} $n,m\ge 3$& $\ssl_{n-2}\oplus\ssl_{m-2}$ & 8 \\
\hline 11  & \begin{picture}(80,66)
    \put(0,51){$\ssp_{2n+2}$}
    \put(40,51){$\ssl_3$}
    \put(10,15){\line(0,1){30}}
    \put(50,15){\line(-4,3){40}}
    \put(50,15){\line(0,1){30}}
    \put(10,15){\circle*{2}}
    \put(10,45){\circle*{2}}
    \put(50,15){\circle*{2}}
    \put(50,45){\circle*{2}}
    \put(0,3){$\ssp_{2n}$}
    \put(40,3){$\ssl_2$}
    \put(52,27){$\pi_1\,\oplus\,1$}
    \end{picture} & $\ssp_{2n-2}$ & $4-\delta_n^0$ \\
\hline 12  & \begin{picture}(100,66)
     \put(20,51){$\ssl_{4}$}
     \put(64,51){$\ssp_{2n+2}$}
     \put(10,15){\line(2,3){20}}
     \put(50,15){\line(-2,3){20}}
     \put(50,15){\line(2,3){20}}
     \put(90,15){\line(-2,3){20}}
     \put(10,15){\circle*{2}}
     \put(50,15){\circle*{2}}
     \put(90,15){\circle*{2}}
     \put(30,45){\circle*{2}}
     \put(70,45){\circle*{2}}
     \put(0,3){$\ssl_2$}
     \put(40,3){$\ssl_2$}
     \put(84,3){$\ssp_{2n}$}
     \end{picture} & $\ssp_{2n-2}$ & $5-\delta_n^0$ \\
\hline 13 & \begin{picture}(120,66)
   \put(15,51){$\ssp_{2n+2}$}
   \put(60,51){$\ssl_{4}$}
   \put(90,51){$\ssp_{2m+2}$}
   \put(7,15){\line(1,2){15}}
   \put(37,15){\line(-1,2){15}}
   \put(37,15){\line(3,4){22}}
   \put(60,15){\line(0,1){30}}
   \put(82,15){\line(-3,4){22}}
   \put(82,15){\line(1,2){15}}
   \put(112,15){\line(-1,2){15}}
   \put(7,15){\circle*{2}}
   \put(37,15){\circle*{2}}
   \put(37,15){\circle*{2}}
   \put(60,15){\circle*{2}}
   \put(82,15){\circle*{2}}
   \put(112,15){\circle*{2}}
   \put(22,45){\circle*{2}}
   \put(60,45){\circle*{2}}
   \put(97,45){\circle*{2}}
   \put(0,3){$\ssp_{2n}$}
   \put(30,3){$\ssl_{2}$}
   \put(60,3){$\cc$}
   \put(75,3){$\ssl_2$}
   \put(105,3){$\ssp_{2m}$}
   \end{picture} & $\ssp_{2n-2}\oplus\ssp_{2m-2}$ &
$7-\delta_n^0-\delta_m^0$\\
\hline 14  & \begin{picture}(115,66)
     \put(0,51){$\ssl_{n+2}$}
     \put(30,51){$\ssp_{2k+2}$}
     \put(90,51){$\ssp_{2m+2}$}
     \put(7,15){\line(0,1){30}}
     \put(37,15){\line(0,1){30}}
     \put(97,15){\line(0,1){30}}
     \put(7,15){\circle*{2}}
     \put(37,15){\circle*{2}}
     \put(97,15){\circle*{2}}
     \put(7,45){\circle*{2}}
     \put(37,45){\circle*{2}}
     \put(97,45){\circle*{2}}
     \put(0,3){$\gl_{n}$}
     \put(30,3){$\ssp_{2k}$}
     \put(90,3){$\ssp_{2m}$}
     \put(67,15){\line(-2,1){60}}
     \put(67,15){\line(-1,1){30}}
     \put(67,15){\line(1,1){30}}
     \put(67,15){\circle*{2}}
     \put(60,3){$\ssl_{2}$}
     \end{picture} $n \ge 1$ & $\gl_{n-2}\oplus\ssp_{2k-2}\oplus\ssp_{2m-2}$
& $7-\delta_n^1-\delta_k^0-\delta_m^0$ \\ \hline 15  &
\begin{picture}(115,66)
     \put(0,51){$\ssl_{n+2}$}
     \put(30,51){$\ssp_{2k+2}$}
     \put(90,51){$\ssp_{2m+2}$}
     \put(7,15){\line(0,1){30}}
     \put(37,15){\line(0,1){30}}
     \put(97,15){\line(0,1){30}}
     \put(7,15){\circle*{2}}
     \put(37,15){\circle*{2}}
     \put(97,15){\circle*{2}}
     \put(7,45){\circle*{2}}
     \put(37,45){\circle*{2}}
     \put(97,45){\circle*{2}}
     \put(0,3){$\ssl_{n}$}
     \put(30,3){$\ssp_{2k}$}
     \put(90,3){$\ssp_{2m}$}
     \put(67,15){\line(-2,1){60}}
     \put(67,15){\line(-1,1){30}}
     \put(67,15){\line(1,1){30}}
     \put(67,15){\circle*{2}}
     \put(60,3){$\ssl_{2}$}
     \end{picture} $n \ge 3$ &
$\ssl_{n-2}\oplus\ssp_{2k-2}\oplus\ssp_{2m-2}$
&$8-\delta_k^0-\delta_m^0  $ \\
\hline
\end{tabular}
}

\pagebreak

\centerline{\large Table 4-3}

\vspace{0.7cm}

\centerline{
\begin{tabular}{|c|c|c|c|}
\hline
 & $\h\subset \g$ & $\s$ &rk\\
\hline 16  & \begin{picture}(134,66)
     \put(37,51){$\ssl_{n+2}$}
     \put(67,51){$\ssp_4$}
     \put(97,51){$\ssp_{2m+2}$}
     \put(30,15){\line(1,2){15}}
     \put(60,15){\line(-1,2){15}}
     \put(60,15){\line(1,2){15}}
     \put(90,15){\line(-1,2){15}}
     \put(90,15){\line(1,2){15}}
     \put(120,15){\line(-1,2){15}}
     \put(30,15){\circle*{2}}
     \put(60,15){\circle*{2}}
     \put(90,15){\circle*{2}}
     \put(120,15){\circle*{2}}
     \put(45,45){\circle*{2}}
     \put(75,45){\circle*{2}}
     \put(105,45){\circle*{2}}
     \put(22,3){$\gl_n$}
     \put(52,3){$\ssl_2$}
     \put(82,3){$\ssl_2$}
     \put(112,3){$\ssp_{2m}$}
     \end{picture} $n\ge 1$ & $\gl_{n-2}\oplus\ssp_{2m-2}$
& $7-\delta_n^1-\delta_m^0$ \\
\hline 17  & \begin{picture}(135,66)
     \put(37,51){$\ssl_{n+2}$}
     \put(67,51){$\ssp_4$}
     \put(97,51){$\ssp_{2m+2}$}
     \put(30,15){\line(1,2){15}}
     \put(60,15){\line(-1,2){15}}
     \put(60,15){\line(1,2){15}}
     \put(90,15){\line(-1,2){15}}
     \put(90,15){\line(1,2){15}}
     \put(120,15){\line(-1,2){15}}
     \put(30,15){\circle*{2}}
     \put(60,15){\circle*{2}}
     \put(90,15){\circle*{2}}
     \put(120,15){\circle*{2}}
     \put(45,45){\circle*{2}}
     \put(75,45){\circle*{2}}
     \put(105,45){\circle*{2}}
     \put(22,3){$\ssl_n$}
     \put(52,3){$\ssl_2$}
     \put(82,3){$\ssl_2$}
     \put(112,1){$\ssp_{2m}$}
     \end{picture} $n\ge 3$ & $\ssl_{n-2}\oplus\ssp_{2m-2}$ & $8-\delta_m^0$
\\
\hline 18 & \begin{picture}(120,66)
   \put(30,51){$\ssp_{2k+4}$}
   \put(60,51){$\ssp_{2n+2}$}
   \put(7,15){\line(1,1){30}}
   \put(37,15){\line(0,1){30}}
   \put(67,15){\line(-1,1){30}}
   \put(67,15){\line(0,1){30}}
   \put(97,15){\line(-1,1){30}}
   \put(7,15){\circle*{2}}
   \put(37,15){\circle*{2}}
   \put(67,15){\circle*{2}}
   \put(97,15){\circle*{2}}
   \put(37,45){\circle*{2}}
   \put(67,45){\circle*{2}}
   \put(0,3){$\ssp_{2k}$}
   \put(30,3){$\ssl_2$}
   \put(60,3){$\ssl_2$}
   \put(90,3){$\ssp_{2n}$}
   \end{picture} $k\ge 1$ & $\ssp_{2k-4}\oplus\ssp_{2n-2}$
& $6-\delta_n^0-\delta_k^1$ \\
\hline 19  & \begin{picture}(135,66)
     \put(0,51){$\ssp_{2n+2}$}
     \put(30,51){$\ssp_{2m+2}$}
     \put(90,51){$\ssp_{2l+2}$}
     \put(120,51){$\ssp_{2k+2}$}
     \put(7,15){\line(0,1){30}}
     \put(37,15){\line(0,1){30}}
     \put(97,15){\line(0,1){30}}
     \put(127,15){\line(0,1){30}}
     \put(7,15){\circle*{2}}
     \put(37,15){\circle*{2}}
     \put(97,15){\circle*{2}}
     \put(127,15){\circle*{2}}
     \put(7,45){\circle*{2}}
     \put(37,45){\circle*{2}}
     \put(97,45){\circle*{2}}
     \put(127,45){\circle*{2}}
     \put(0,3){$\ssp_{2n}$}
     \put(30,3){$\ssp_{2m}$}
     \put(90,3){$\ssp_{2l}$}
     \put(120,3){$\ssp_{2k}$}
     \put(67,15){\line(-2,1){60}}
     \put(67,15){\line(-1,1){30}}
     \put(67,15){\line(1,1){30}}
     \put(67,15){\line(2,1){60}}
     \put(67,15){\circle*{2}}
     \put(60,3){$\ssl_2$}
     \end{picture}
&$\ssp_{2n-2}\oplus\ssp_{2m-2}\oplus$
& $8-\delta_n^0-\delta_m^0-$ \\
&&$\oplus\ssp_{2l-2}\oplus\ssp_{2k-2}$&$-\delta_l^0-\delta_k^0$\\
\hline 20  & \begin{picture}(135,66)
     \put(0,51){$\ssp_{2n+2}$}
     \put(30,51){$\ssp_{2k+2}$}
     \put(90,51){$\ssp_{4}$}
     \put(120,51){$\ssp_{2m+2}$}
     \put(7,15){\line(0,1){30}}
     \put(37,15){\line(0,1){30}}
     \put(97,15){\line(0,1){30}}
     \put(127,15){\line(0,1){30}}
     \put(7,15){\circle*{2}}
     \put(37,15){\circle*{2}}
     \put(97,15){\circle*{2}}
     \put(127,15){\circle*{2}}
     \put(7,45){\circle*{2}}
     \put(37,45){\circle*{2}}
     \put(97,45){\circle*{2}}
     \put(127,45){\circle*{2}}
     \put(0,3){$\ssp_{2n}$}
     \put(30,3){$\ssp_{2k}$}
     \put(90,3){$\ssl_2$}
     \put(120,3){$\ssp_{2m}$}
     \put(67,15){\line(-2,1){60}}
     \put(67,15){\line(-1,1){30}}
     \put(67,15){\line(1,1){30}}
     \put(97,15){\line(1,1){30}}
     \put(67,15){\circle*{2}}
     \put(60,3){$\ssl_2$}
     \end{picture}
&$\ssp_{2n-2}\oplus\ssp_{2m-2}\oplus\ssp_{2k-2}$ &
$8-\delta_n^0-\delta_m^0-\delta_k^0$ \\
\hline 21  & \begin{picture}(135,66)
     \put(15,51){$\ssp_{2n+2}$}
     \put(45,51){$\ssp_4$}
     \put(75,51){$\ssp_4$}
     \put(105,51){$\ssp_{2m+2}$}
     \put(7,15){\line(1,2){15}}
     \put(37,15){\line(-1,2){15}}
     \put(37,15){\line(1,2){15}}
     \put(67,15){\line(-1,2){15}}
     \put(67,15){\line(1,2){15}}
     \put(97,15){\line(-1,2){15}}
     \put(97,15){\line(1,2){15}}
     \put(127,15){\line(-1,2){15}}
     \put(7,15){\circle*{2}}
     \put(37,15){\circle*{2}}
     \put(67,15){\circle*{2}}
     \put(97,15){\circle*{2}}
     \put(127,15){\circle*{2}}
     \put(22,45){\circle*{2}}
     \put(52,45){\circle*{2}}
     \put(82,45){\circle*{2}}
     \put(112,45){\circle*{2}}
     \put(0,3){$\ssp_{2n}$}
     \put(30,3){$\ssl_2$}
     \put(60,3){$\ssl_2$}
     \put(90,3){$\ssl_2$}
     \put(120,3){$\ssp_{2m}$}
     \end{picture} & $\ssp_{2n-2}\oplus\ssp_{2m-2}$ &
$8-\delta_n^0-\delta_m^0$ \\
\hline 22  & \begin{picture}(100,66)
     \put(20,51){$\ssp_{4}$}
     \put(64,51){$\ssp_{2n+2}$}
     \put(10,15){\line(2,3){20}}
     \put(50,15){\line(-2,3){20}}
     \put(50,15){\line(2,3){20}}
     \put(90,15){\line(-2,3){20}}
     \put(10,15){\circle*{2}}
     \put(50,15){\circle*{2}}
     \put(90,15){\circle*{2}}
     \put(30,45){\circle*{2}}
     \put(70,45){\circle*{2}}
     \put(0,3){$\cc$}
     \put(40,3){$\ssl_2$}
     \put(84,3){$\ssp_{2n}$}
     \put(36,36){{\tiny $\pi_1\oplus 1\oplus 1$}}
     \end{picture} & $\ssp_{2n-2}$ & $4-\delta_n^0$ \\
\hline
\end{tabular}
} \pagebreak

\centerline{\large Table 4-4}

\vspace{0.7cm}

\centerline{
\begin{tabular}{|c|c|c|c|}
\hline
 & $\h\subset \g$ & $\s$ &rk\\
\hline 23  & \begin{picture}(125,66)
     \put(40,51){$\ssp_6$}
     \put(84,51){$\ssp_{2n+2}$}
     \put(4,51){$\ssp_4$}
     \put(30,15){\line(-2,3){20}}
     \put(30,15){\line(2,3){20}}
     \put(70,15){\line(-2,3){20}}
     \put(70,15){\line(2,3){20}}
     \put(110,15){\line(-2,3){20}}
     \put(10,45){\circle*{2}}
     \put(30,15){\circle*{2}}
     \put(70,15){\circle*{2}}
     \put(110,15){\circle*{2}}
     \put(50,45){\circle*{2}}
     \put(90,45){\circle*{2}}
     \put(20,3){$\ssp_{4}$}
     \put(60,3){$\ssl_2$}
     \put(104,3){$\ssp_{2n}$}
     \end{picture} & $\ssp_{2n-2}$ & $7-\delta_n^0$ \\
\hline 24  & \begin{picture}(80,66)
     \put(0,51){$G_2$}
     \put(40,51){$\ssl_3$}
     \put(30,15){\line(-2,3){20}}
     \put(30,15){\line(2,3){20}}
     \put(30,15){\circle*{3}}
     \put(10,45){\circle*{3}}
     \put(50,45){\circle*{3}}
     \put(20,3){$\ssl_3$}
     \end{picture}& 0&4\\
\hline 25  & \begin{picture}(80,66)
    \put(0,51){$\ssl_{n+3}$}
    \put(40,51){$\ssl_3$}
    \put(10,15){\line(0,1){30}}
    \put(50,15){\line(-4,3){40}}
    \put(50,15){\line(0,1){30}}
    \put(10,15){\circle*{2}}
    \put(10,45){\circle*{2}}
    \put(50,15){\circle*{2}}
    \put(50,45){\circle*{2}}
    \put(0,3){$\gl_n$}
    \put(40,3){$\ssl_3$}
    \end{picture} $n\ge 2$ & $\gl_{n-3}$ & $7-\delta_n^2$ \\
\hline 26  & \begin{picture}(80,66)
    \put(0,51){$\ssl_{n+3}$}
    \put(40,51){$\ssl_3$}
    \put(10,15){\line(0,1){30}}
    \put(50,15){\line(-4,3){40}}
    \put(50,15){\line(0,1){30}}
    \put(10,15){\circle*{2}}
    \put(10,45){\circle*{2}}
    \put(50,15){\circle*{2}}
    \put(50,45){\circle*{2}}
    \put(0,3){$\ssl_n$}
    \put(40,3){$\ssl_3$}
    \end{picture} $n\ge 4$ & $\ssl_{n-3}$ & 8 \\
\hline 27  & \begin{picture}(80,66)
     \put(0,51){$\ssl_{n+1}$}
     \put(40,51){$\ssl_{n}$}
     \put(30,15){\line(-2,3){20}}
     \put(30,15){\line(2,3){20}}
     \put(30,15){\circle*{2}}
     \put(10,45){\circle*{2}}
     \put(50,45){\circle*{2}}
     \put(20,3){$\ssl_n$}
     \end{picture} $n\ge 2$ & 0 & $2n-1$ \\
\hline 28  & \begin{picture}(80,66)
    \put(0,51){$\ssp_8$}
    \put(40,51){$\ssp_6$}
    \put(10,15){\line(0,1){30}}
    \put(50,15){\line(-4,3){40}}
    \put(50,15){\line(0,1){30}}
    \put(10,15){\circle*{2}}
    \put(10,45){\circle*{2}}
    \put(50,15){\circle*{2}}
    \put(50,45){\circle*{2}}
    \put(0,3){$\ssl_2$}
    \put(40,3){$\ssp_6$}
    \end{picture} & 0 & 7 \\
\hline 29  & \begin{picture}(80,66)
     \put(0,51){$\so_7$}
     \put(40,51){$G_2$}
     \put(30,15){\line(-2,3){20}}
     \put(30,15){\line(2,3){20}}
     \put(30,15){\circle*{2}}
     \put(10,45){\circle*{2}}
     \put(50,45){\circle*{2}}
     \put(20,3){$G_2$}
     \end{picture}&0&5\\
\hline 30  & \begin{picture}(120,66)
    \put(0,51){$\ssl_3$}
    \put(40,51){$\ssl_3$}
    \put(80,51){$\ssl_3$}
    \put(50,15){\line(4,3){40}}
    \put(50,15){\line(0,1){30}}
    \put(50,15){\line(-4,3){40}}
    \put(50,15){\circle*{2}}
    \put(10,45){\circle*{2}}
    \put(50,45){\circle*{2}}
    \put(90,45){\circle*{2}}
    \put(40,3){$\ssl_3$}
    \end{picture}&0&6\\
\hline
\end{tabular}
} \pagebreak

\normalsize

\bigskip

{\hspace{60mm} {\small \textsf Translated by O. V. Chuvashova}}

\bigskip
\bigskip

\end{document}